\documentclass{amsart}
\usepackage{amssymb}
\usepackage{amsthm}
 
\newtheorem{prop}{Proposition}
\theoremstyle{definition} \newtheorem{defin}{Definition}[section]
\newtheorem{ex}{Example}[section]
\theoremstyle{remark} \newtheorem{rem}{Remark}[section]
\newcommand{\pn}{\par\noindent} \newcommand{\pmn}{\par\medskip\noindent}

\begin{document}
\title{Real Morse polynomials of degrees 5 and 6}
\author{Yury Kochetkov}
\date{}
\begin{abstract} A real polynomial $p$ of degree $n$ is called a
Morse polynomial if its derivative has $n-1$ pairwise different
real roots and values of $p$ in these roots (critical values) are also
pairwise different. The plot of such polynomial is called a
"snake". By enumerating  critical points and critical values in
the increasing order we construct a permutation
$a_1,\ldots,a_{n-1}$, where $a_i$ is the number of polynomial's
value in $i$-th critical point. This permutation is called the
\emph{passport} of the snake (polynomial). In this work for Morse
polynomials of degrees 5 and 6 we describe the partition of the
coefficient space into domains of constant passport.
\end{abstract} \email{yukochetkov@hse.ru, yuyukochetkov@gmail.com}
\maketitle
\section{Introduction}

\begin{defin} A real polynomial will be called a real Morse
polynomial if:
\begin{itemize} \item all its critical points (i.e. roots of its
derivative) are real and have multiplicity 2 (i.e. all roots of
derivative are real and simple); \item all its critical values are
pairwise different; \item the leading coefficient is 1.
\end{itemize}
\end{defin}\pn The plot of
such polynomial will be called (after Arnold \cite{Ar}) a "snake".
A snake of order $n$ is the plot of a  real Morse polynomial of
degree $n+1$. It has $n$ alternating local minima and maxima. \pmn
Let's enumerate critical points in the increasing order and
critical values also in the increasing order. The number of the
first critical value will be called the \emph{level} of a snake.
To an $n$-snake we will correspond a permutation $a_1,\ldots,a_n$,
where $a_i$ is the number of critical value in the $i$-th critical
point. Thus defined permutation will be called the \emph{passport}
of a snake (polynomial). Two snakes with the same passports will
be considered equal and two real Morse polynomials with equal
passports will be called $s$-equivalent.

\begin{defin} A permutation $a_1,\ldots,a_n$ is called \emph{alternating},
if a short sequence $a_i,a_{i+1},a_{i+2}$ is not monotone for any $i$,
$1\leqslant i\leqslant n-2$. \end{defin}

\pn The passport of a snake is an alternating permutation. The
level of an alternating permutation is its first element. As we
study plots of normalized polynomials, then the level of a snake
of even order $n$ is greater, than 1 (because the first extremum
of a normalized polynomial of odd degree is maximum), and the
level of a snake of odd order $n$ is less, than $n$ (because the
first extremum of a normalized polynomial of even order is
minimum). Thus, we will study not all alternating permutations,
but \emph{proper} alternating permutations. The level of such
permutation from $S_n$ is greater, than 1, for even $n$, and less,
than $n$, for odd $n$. \begin{rem} The answer to the natural
question: is it true that any proper alternating permutation is a
passport of some snake?, is positive (see \cite{IZ} and
bibliography there).\end{rem} \pn In what follows PAP will be a
proper alternating permutation.

\begin{ex} There are 5 PAPs of order 4:
$(4,2,3,1)$, $(3,2,4,1)$, $(4,1,3,2)$, $(3,1,4,2)$ and $(2,1,4,3)$,
and there are 5 snakes of order 4:
\[\begin{picture}(380,65) \multiput(0,14)(10,0){38}{\line(1,0){6}}
\multiput(0,28)(10,0){38}{\line(1,0){6}}
\multiput(0,42)(10,0){38}{\line(1,0){6}}
\multiput(0,56)(10,0){38}{\line(1,0){6}}
\qbezier(0,5)(5,56)(10,56) \qbezier(10,56)(13,56)(17,42)
\qbezier(17,42)(21,28)(24,28) \qbezier(24,28)(27,28)(30,35)
\qbezier(30,35)(33,42)(36,42) \qbezier(36,42)(39,42)(43,28)
\qbezier(43,28)(47,14)(50,14) \qbezier(50,14)(55,14)(60,62)

\put(30,2){\scriptsize 1}

\qbezier(80,5)(85,42)(90,42) \qbezier(90,42)(93,42)(96,35)
\qbezier(96,36)(99,28)(102,28) \qbezier(102,28)(105,28)(109,42)
\qbezier(109,42)(113,56)(116,56) \qbezier(116,56)(120,56)(124,35)
\qbezier(124,35)(128,14)(132,14) \qbezier(132,14)(137,14)(142,62)

\put(110,2){\scriptsize 2}

\qbezier(160,5)(165,56)(170,56) \qbezier(170,56)(174,56)(178,35)
\qbezier(178,35)(182,14)(186,14) \qbezier(186,14)(189,14)(193,28)
\qbezier(193,28)(197,42)(200,42) \qbezier(200,42)(203,42)(206,35)
\qbezier(206,35)(209,28)(212,28) \qbezier(212,28)(216,28)(220,62)

\put(190,2){\scriptsize 3}

\qbezier(240,5)(245,42)(250,42) \qbezier(250,42)(253,42)(257,28)
\qbezier(257,28)(261,14)(264,14) \qbezier(264,14)(268,14)(272,35)
\qbezier(272,35)(276,56)(280,56) \qbezier(280,56)(283,56)(287,42)
\qbezier(287,42)(291,28)(294,28) \qbezier(294,28)(298,28)(302,62)

\put(272,2){\scriptsize 4}

\qbezier(320,5)(324,28)(328,28) \qbezier(328,28)(331,28)(334,21)
\qbezier(334,21)(337,14)(340,14) \qbezier(340,14)(344,14)(348,35)
\qbezier(348,35)(352,56)(356,56) \qbezier(356,56)(359,56)(362,49)
\qbezier(362,49)(365,42)(368,42) \qbezier(368,42)(372,42)(376,62)

\put(348,2){\scriptsize 5}
\end{picture}\]
\begin{center}{\Large Figure 1}\end{center}
\end{ex}

\begin{rem} There are 16 PAPs of order 5 and 16 snakes
of order 5. \end{rem}

\section{The enumeration of PAPs}
\pn Let us introduce a procedure, that transform a PAP of order
$n$ and level $m$ into PAP of order $n+1$: we add the first
element $k$ and all elements $<k$ will be the same and all
elements $\geqslant k$ will be increased by 1. If $n$ is even then
$1\leqslant k\leqslant m$. If $n$ is odd then $m<k\leqslant n+1$.
For example,
$$(3,1,5,2,6,4)\quad\text{add 2}\quad\Rightarrow (2,4,1,6,3,7,5),$$
$$(2,4,1,5,3)\quad\text{add 4}\quad\Rightarrow (4,2,5,1,6,3).$$
Obviously, this procedure transforms PAP into PAP.\pmn The inverse
procedure --- the deletion of the first element: we delete the
first element $m$, all elements $<m$ remain the same, all elements
$>m$ are decreased by 1. \pmn Thus, any PAP of order $n+1$ can be
uniquely obtained from some PAP of order $n$. \pmn Let $s(n,m)$ be
the number of PAP of order $n$ and level $m$. If $n$ is even then
$$s(n+1,m)=s(n,m)+\ldots+s(n,n)\,.\eqno(1)$$ If $n$ is odd then
$$s(n+1,m)=s(n,1)+\ldots+s(n,m-1)\,.\eqno(2)$$ Thus,
$$\begin{array}{l} s(6,2)=s(5,1)\\ s(6,3)=s(5,1)+s(5,2)\\ s(6,4)=
s(5,1)+s(5,2)+s(5,3)\\ s(6,5)=s(5,1)+s(5,2)+s(5,3)+s(5,4)\\
s(6,6)=s(5,1)+s(5,2)+s(5,3)+s(5,4)\end{array}$$ (we remind that
$s(5,5)=0$ and $s(6,1)=0$). Analogously,
$$\begin{array}{l}s(5,1)=s(4,2)+s(4,3)+s(4,4)\\ s(5,2)=s(4,2)+s(4,3)+
s(4,4)\\ s(5,3)=s(4,3)+s(4,4)\\ s(5,4)=s(4,4) \end{array}$$ (we
remind that $s(4,1)=0$ and $s(5,5)=0$).  \pmn Now we can construct
the triangle of PAPs --- the Euler-Bernoulli triangle (see
\cite{Ar}). In $n$-th row we write numbers $s(n,1),\ldots,s(n,n)$
from left to right. The case $s(1,1)$ is special, we put
$s(1,1)=1$. From (1) and (2) we have that an element in an even
row is the sum of elements in the previous row that are to the
left of our element, and an element in an odd row is the sum of
elements in the previous row that are to the right of our element.
Below are presented the first 5 rows of our triangle:
$$\begin{array}{ccccccccc}&&&&1&&&&\\ &&&0&&1&&&\\ &&1&&1&&0&&\\
&0&&1&&2&&2&\\ 5&&5&&4&&2&&0\end{array}$$

\section{Polynomials of degree 5}
\pn Our aim is to describe the partition of the space of real
Morse polynomials of degrees 5 and 6 into components of
s-equivalency. In this section we will consider polynomials of
degree 5. \pmn At first let us introduce a convenient
parametrization. Let $p$ be a real Morse polynomial of degree 5.
We will assume that: a) roots of its derivative are non positive
and one root is zero; b) the sum of roots of $p'$ is $-3$. It
means that $p'=x^4+3x^3+bx^2+cx$, where $b>0$ and $c>0$. \pmn Let
us note that if coefficients of the polynomial $q=x^3+3x^2+bx+c$
are positive, then its roots are negative (and their sum is $-3$).
Indeed, let $q=(x+x_1)(x+x_2)(x+x_3)$, then either all $x_i>0$, or
$x_1>0$ and $x_2<0,x_3<0$. Let us consider the second case and let
$y_2=-x_2>0$ and $y_3=-x_3>0$, then $x_1=3+y_2+y_3$ and
$y_2y_3>x_1(y_2+y_3)=3(y_2+y_3)+y_2^2+2y_2y_3+y_3^2$.
Contradiction. \pmn We will work with the polynomial
$q=x^3+3x^2+bx+c$ with negative pairwise different roots. Then
$b\leqslant 3$ and $c\leqslant 1$. The discriminant $dq$ of $q$ is
$dq=54bc-108c-4b^3+9b^2-27c^2$ and the plot of the curve $dq=0$ in
the rectangle $[0,3]\times [0,1]$ is presented below:
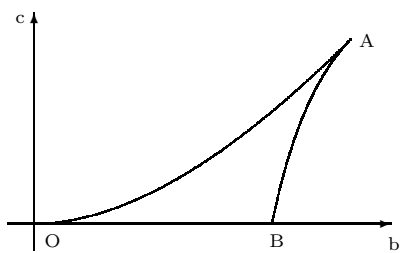
\[\begin{picture}(200,100) \put(0,15){\vector(1,0){145}}
\put(10,5){\vector(0,1){90}} \put(144,5){\scriptsize b}
\put(3,91){\scriptsize c} \qbezier(10,15)(60,15)(130,85)
\qbezier(100,15)(110,65)(130,85) \put(133,82){\scriptsize A}
\put(99,6){\scriptsize B} \put(14,6){\scriptsize O}
\end{picture}\] The point $A$ has coordinates $(3,1)$ and is a cusp of
the curve $dq=0$. The line $y=x-2$ is the tangent line to the
curve at this point. The point $B$ has coordinates $(9/4,0)$ and
the line $y=\frac 32 x-\frac{27}{8}$ is the tangent line to the
curve at this point. In the origin the curve $dq=0$ is tangent to
axis $OX$. For points in the curvilinear triangle $OAB$ the
polynomial $q$ has three real roots. \pmn We must also demand that
critical values of polynomial $p=\int xq(x)\,dx$ are pairwise
different. Let us consider two conditions
$$g(b,c)=128b^3-1998b^2-216c^2+1512bc+3645b-729c\neq 0$$  and
$$h(b,c)=640b^3-1350b^2+5832c^2-9720bc+18225c\neq 0\,.$$ The
satisfaction of the first condition guarantee that values of the
polynomial $p$ in roots of $q$ are pairwise different and the
satisfaction of the second condition guarantee that these values
are nonzero. \pmn The curve $g=0$ is an arc inside $OAB$ that
connect $A$ with the point $D$ with coordinates $(135/64,0)$.
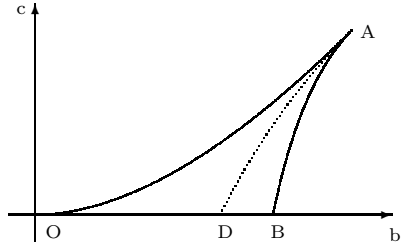
\[\begin{picture}(200,100) \put(0,15){\vector(1,0){145}}
\put(10,5){\vector(0,1){90}} \put(144,5){\scriptsize b}
\put(3,91){\scriptsize c} \qbezier(10,15)(60,15)(130,85)
\qbezier(100,15)(110,65)(130,85) \put(133,82){\scriptsize A}
\put(99,6){\scriptsize B} \put(14,6){\scriptsize O}
\qbezier[50](130,85)(100,55)(80,15) \put(79,6){\scriptsize D}
\end{picture}\] The line $y=\frac 54 x-\frac{675}{256}$ is tangent to
the curve $g=0$ at the point $D$. \pmn The plot of the curve $h=0$
inside $OAB$ is presented below:
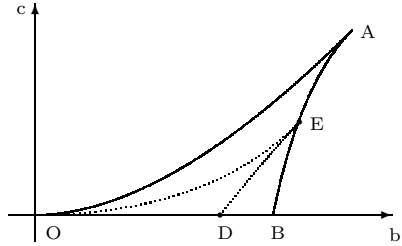
\[\begin{picture}(200,100) \put(0,15){\vector(1,0){145}}
\put(10,5){\vector(0,1){90}} \put(144,5){\scriptsize b}
\put(3,91){\scriptsize c} \qbezier(10,15)(60,15)(130,85)
\qbezier(100,15)(110,65)(130,85) \put(133,82){\scriptsize A}
\put(99,6){\scriptsize B} \put(14,6){\scriptsize O}
\put(79,6){\scriptsize D} \put(114,47){\scriptsize E}
\qbezier[40](110,50)(100,40)(80,15) \put(110,50){\circle*{2}}
\qbezier[60](110,50)(75,15)(10,15) \put(80,15){\circle*{2}}
\end{picture}\] The point $E$ with coordinates $(45/16,25/32)$ is
the cusp point of the curve $h=0$. Lines $y=\frac 56
x-\frac{25}{16}$ and $y= \frac 54 x-\frac{175}{64}$ are tangent
lines to curves $h=0$ and $dq=0$ at this point. Lines $g=0$ and
$h=0$ are tangent to each other at the point $D$. \pmn In figure
below we demonstrate how the curvilinear triangle $OAB$ is
partitioned into five curvilinear triangles by curves $g=0$ and
$h=0$.
\[\begin{picture}(180,135) \put(10,10){\line(1,0){110}}
\put(90,10){\line(2,3){80}} \qbezier(10,10)(90,10)(170,130)
\put(120,10){\line(1,3){20}} \qbezier(140,70)(150,100)(170,130)
\qbezier(140,70)(110,40)(90,10) \qbezier(140,70)(80,10)(10,10)
\put(8,1){\scriptsize O} \put(88,1){\scriptsize D}
\put(118,1){\scriptsize B} \put(173,127){\scriptsize A}
\put(143,65){\scriptsize E} \put(110,48){\scriptsize F}
\end{picture}\]
\begin{center}{\Large Figure 2}\end{center}
\pmn
The point $F$ --- the intersection point of curves
$g=0$ and $h=0$ has coordinates $\approx(2.73,0.72)$. \pmn All five
snakes of order 4 are presented at Figure 1. Let
$$p(x)=5\int (x^4+3x^3+bx^2+cx)\,dx=x^5+\frac{15}{4}\,x^4+
\frac{5b}{3}\,x^3+\frac{5c}{2}\,x^2\,.$$
\begin{itemize} \item If a point with coordinates $(b,c)$ is in triangle
$OAF$, then the plot of $p$ is the first snake; \item if a point
with coordinates $(b,c)$ is in triangle $AEF$, then the plot of
$p$ is the second snake; \item if a point with coordinates $(b,c)$
is in triangle $ODF$, then the plot of $p$ is the third snake;
\item if a point with coordinates $(b,c)$ is in triangle $DEF$,
then the plot of $p$ is the forth snake; \item if a point with
coordinates $(b,c)$ is in triangle $BDE$, then the plot of $p$ is
the fifth snake.\end{itemize}
\begin{rem} If a real polynomial of degree 5 has four real critical
points but only three critical values, then the plot of this
polynomial is a "degenerate"{} snake. On the figure below are
presented all five degenerate snakes of order 4.
\[\begin{picture}(390,50) \multiput(0,14)(10,0){39}{\line(1,0){6}}
\multiput(0,28)(10,0){39}{\line(1,0){6}}
\multiput(0,42)(10,0){39}{\line(1,0){6}}
\qbezier(0,5)(7,42)(14,42) \qbezier(14,42)(17,42)(21,28)
\qbezier(21,28)(25,14)(28,14) \qbezier(28,14)(31,14)(35,28)
\qbezier(35,28)(39,42)(42,42) \qbezier(42,42)(45,42)(48,35)
\qbezier(48,35)(51,28)(54,28) \qbezier(54,28)(57,28)(60,48)

\put(30,2){\scriptsize 1}

\qbezier(80,5)(87,42)(94,42) \qbezier(94,42)(97,42)(100,35)
\qbezier(100,35)(103,28)(106,28) \qbezier(106,28)(109,28)(112,35)
\qbezier(112,35)(115,42)(118,42) \qbezier(118,42)(121,42)(125,28)
\qbezier(125,28)(129,14)(132,14) \qbezier(132,14)(139,14)(146,48)

\put(110,2){\scriptsize 2}

\qbezier(166,5)(173,42)(180,42) \qbezier(180,42)(183,42)(187,28)
\qbezier(187,28)(191,14)(194,14) \qbezier(194,14)(197,14)(200,21)
\qbezier(200,21)(203,28)(206,28) \qbezier(206,28)(209,28)(212,21)
\qbezier(212,21)(215,14)(218,14) \qbezier(218,14)(225,14)(232,48)

\put(196,2){\scriptsize 3}

\qbezier(252,5)(256,28)(260,28) \qbezier(260,28)(263,28)(266,21)
\qbezier(266,21)(269,14)(272,14) \qbezier(272,14)(275,14)(279,28)
\qbezier(279,28)(283,42)(286,42) \qbezier(286,42)(289,42)(293,28)
\qbezier(293,28)(297,14)(300,14) \qbezier(300,14)(304,14)(310,48)

\put(282,2){\scriptsize 4}

\qbezier(330,5)(334,28)(338,28) \qbezier(338,28)(341,28)(344,21)
\qbezier(344,21)(347,14)(350,14) \qbezier(350,14)(353,14)(357,28)
\qbezier(357,28)(361,42)(364,42) \qbezier(364,42)(367,42)(370,35)
\qbezier(370,35)(373,28)(376,28) \qbezier(376,28)(380,28)(385,48)

\put(360,2){\scriptsize 5} \end{picture}\]

\begin{center}{\Large Figure 3}\end{center}
\pmn Arcs, that separate curvilinear triangles in Figure 2,
correspond to degenerate snakes:
\begin{itemize} \item the arc $AF$ --- to the second degenerate snake;
\item the arc $OF$ --- to the third degenerate snake; \item the
arc $DF$ --- to the first degenerate snake; \item the arc $EF$
--- to the forth degenerate snake; \item the arc $DE$ --- to the fifth
degenerate snake. \end{itemize} \end{rem}

\section{Polynomials of degree 6}
\pn We will introduce the analogous parametrization. Let $p$ be a
real Morse polynomial. We will assume, that: a) all roots of its
derivative are non positive and pairwise different; b) one root is
zero; c) the sum of roots is $-4$. It means that
$p'=x(x^4+4x^3+ax^2+bx+x)$, where $0<a\leqslant 6$, $0<b\leqslant
4$, $0<c\leqslant 1$.
\begin{prop} If coefficients $a,b,c$ of real polynomial
$q=x^4+4x^3+ax^2+bx+c$ are positive, then its roots are negative.
\end{prop}  \begin{proof} The number of positive roots must be even
(otherwise $c<0$). All roots cannot be positive, because their sum is $-4$.
Hence, we must consider the case of two positive roots $x_1$ and $x_2$
and two negative roots $-y_1$ and $-y_2$. \pmn Let $y_1+y_2=2k$ and
$x_1+x_2=2l$. Then $k-l=2$ and
$$\left\{\begin{array}{l} x_1x_2+y_1y_2>(x_1+x_2)(y_1+y_2)\\
x_1x_2(y_1+y_2)>y_1y_2(x_1+x_2)\end{array}\right.\Rightarrow
\left\{\begin{array}{l} x_1x_2+y_1y_2>4kl\\
kx_1x_2>ly_1y_2\end{array} \right.$$ As $l^2\geqslant x_1x_2$, then
$$\left\{\begin{array}{l} y_1y_2>4kl-l^2=3k^2-4k-4\\y_1y_2<kl=k^2-2k
\end{array}\right.\Rightarrow 3k^2-4k-4<k^2-2k\Rightarrow k^2-k-2<0
\Rightarrow k<2.$$ Contradiction. \end{proof} \pmn The space of
coefficients here will the the parallelepiped $\Pi$ in the axes $a,b,c$:
$\Pi=(0,6]\times (0,4]\times (0,1]$. We will study values $p(x_i)$
of the polynomial $p=\int xq(x)\,dx$, where $q=x^4+4x^3+ax^2+bx+c$ and
$x_1,x_2,x_3,x_4$ are roots of $q$. We restrict the study to those
domain of $\Pi$, where all roots of $q$ are real. \pmn Let
\begin{multline*}d(a,b,c)=16a^4c-4a^3b^2-64a^3c+16a^2b^2-320a^2bc-128a^2c^2
+72ab^3+144ab^2c+\\
+1152abc+2304ac^2-27b^4-256b^3-96b^2c-768bc^2+256c^3-6912c^2\end{multline*}
be the discriminant of $q$. We will consider sections of $\Pi$ by
planes $c=\gamma$ and will study curves $d(a,b,\gamma)=0$ in
rectangle $R=\{0\leqslant a\leqslant 6,0\leqslant b\leqslant 4\}$.
Curves $d(a,b,1)=0$ and $d(a,b,3/4)=0$ in the rectangle
$\{3\leqslant a\leqslant 6, 0\leqslant b\leqslant 4\}$ are
presented below in Figure 4.
\[\begin{picture}(365,170) \put(0,70){\vector(1,0){75}}
\put(5,65){\vector(0,1){55}} \qbezier[50](5,110)(35,110)(65,110)
\qbezier[35](65,70)(65,90)(65,110) \put(25,70){\line(1,1){40}}
\put(70,62) {\scriptsize a} \put(-3,115){\scriptsize b}
\put(23,61){\scriptsize A}

\put(100,70){\vector(1,0){75}} \put(105,65){\vector(0,1){55}}
\qbezier[50](105,110)(135,110)(165,110)
\qbezier[35](165,70)(165,90)(165,110) \put(125,70){\line(4,3){40}}
\put(170,62) {\scriptsize a} \put(97,115){\scriptsize b}
\put(123,61){\scriptsize B} \put(168,98){\scriptsize C}
\put(157,93){$\square$}

\qbezier[100](195,5)(280,5)(365,5)
\qbezier[100](195,165)(280,165)(365,165)
\qbezier[100](195,5)(195,80)(195,165)
\qbezier[100](365,5)(365,80)(365,165)
\qbezier(200,10)(282,61)(365,112) \qbezier(238,5)(300,82)(360,160)
\qbezier(200,10)(290,75)(360,160)  \end{picture}\]
\begin{center}{\Large Figure 4}\end{center}
\pmn  Here point $A$ has coordinates $\approx(3.76,0)$, point $B$
--- $\approx(3.82,0)$, point $C$ --- $\approx(6,3.62)$. The
picture on the right demonstrate a segment of the curve
$d(a,b,3/4)=0$ (the small square on the middle picture) in the
rectangle $\{5.63\leqslant a \leqslant 5.81,3.39\leqslant
b\leqslant 5.81\}$ (the correct scale). In this segment the curve
$d(a,b,3/4)=0$ has two cusps and one self-crossing point. The
curve $d(a,b,3/4)=0$ divides $R$ into three domains: above the
curve --- here $q$ has two real and two complex roots, below the
curve --- here $q$ has only complex roots and inside the
curvilinear triangle --- here $q$ has four real roots. \pmn With
decreasing of $c$ the curvilinear triangle is moving down and to
the left and also is enlarged. In Figures below are presented: a)
plot of the curve $d(a,b,1/16)=0$ in the rectangle $[3\leqslant
a\leqslant 6,0\leqslant b\leqslant 3]$ (the left figure); b) the
plot of the curve $d(a,b,0)=0$ in the rectangle $[0\leqslant
a\leqslant 6,0\leqslant b\leqslant 3]$:
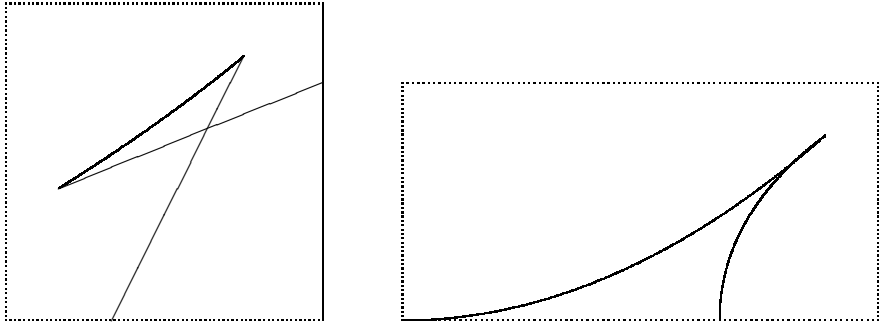
\[\begin{picture}(330,130) \qbezier[70](0,5)(60,5)(120,5)
\qbezier[70](0,125)(60,125)(120,125)
\qbezier[70](0,5)(0,65)(0,125) \qbezier(120,5)(120,65)(120,125)
\put(40,5){\line(1,2){50}} \put(20,55){\line(5,2){100}}
\qbezier(20,55)(60,80)(90,105)

\qbezier[90](150,5)(240,5)(330,5)
\qbezier[90](150,95)(240,95)(330,95)
\qbezier[55](150,5)(150,50)(150,95)
\qbezier[55](330,5)(330,50)(330,95)
\qbezier(270,5)(270,45)(310,75) \qbezier(150,5)(230,5)(310,75)
\end{picture}\]  On the right figure the curve vertically intersects the
axis $OX$ at the point $(4,0)$ and the cusp has coordinates
$\approx (5.33,2.37)$. \pmn In what follows we will study the
interior of the curvilinear triangle only. This triangle will be
called \emph{the main triangle} and will be denoted
$\Delta_\gamma$ (i.e. the main triangle in the plane $c=\gamma$).
\pmn The condition that values of $p$ at roots of $q$ are pairwise
different is of the form $s(a,b,c)\neq 0$, where
$$\begin{array}{l}
s=37500a^7b-15625a^6b^2-688000a^6b-90000a^6c-\\
-480000a^5b^2-450000a^5bc+3624960a^5b +1651200a^5c+\\
+600000a^4b^3+187500a^4b^2c+13363200a^4b^2+13824000a^4bc-\\
-5898240a^4b+1080000a^4c^2-8699904a^4c-125000a^3b^4
-9344000a^3b^3-\\
-8160000a^3b^2c-73662464a^3b^2+1800000a^3bc^2-163184640a^3bc
-28416000a^3c^2+\\+14155776a^3c
+1440000a^2b^4+1200000a^2b^3c+3932160a^2b^3-750000a^2b^2c^2+\\
+126259200a^2b^2c
+116391936a^2b^2-24576000a^2bc^2+723517440a^2bc-4320000a^2c^3+\\
+299630592a^2c^2
+19046400ab^4-46080000ab^3c+154140672ab^3+23040000ab^2c^2-\\
-438829056ab^2c -2400000abc^3-47185920abc^2
-1056964608abc+70656000ac^3-\\
-1264582656ac^2
-4096000b^5+7680000b^4c-66322432b^4-4800000b^3c^2+\\
+96337920b^3c
+1000000b^2c^3 -8601600b^2c^2+276824064b^2c-23552000bc^3+\\
+421527552bc^2
-203423744c^3-268435456b^3+1811939328c^2+5760000c^4.
\end{array}$$  The intersection of the main triangle $\Delta_\gamma$
and the surface $s=0$ is of the form
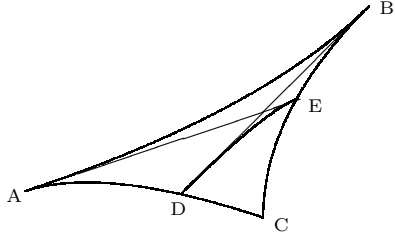
\[\begin{picture}(150,100) \put(10,20){\line(3,1){104}}
\put(140,90){\line(-1,-1){71}} \qbezier(10,20)(100,50)(140,90)
\qbezier(140,90)(100,50)(100,10) \qbezier(10,20)(40,30)(100,10)
\qbezier(69,19)(100,50)(113,55) \put(3,16){\scriptsize A}
\put(144,87){\scriptsize B} \put(104,5){\scriptsize C}
\put(65,11){\scriptsize D} \put(117,50){\scriptsize E}
\end{picture}\] Here $ABC$ is the main triangle, $D$ and $E$ are
cusps of the curve $s(a,b,\gamma)=0$ (both are on sides of the
main triangle), arcs of curve $s(a,b,\gamma)=0$ are outgoing  from
cusp points $A$ and $B$ to the interior of the main triangle.
\pmn We will also need the condition that values of $p$ at roots of
the polynomial $q$ are nonzero. This condition is of the form
$z(a,b,c)\neq 0$, where
\begin{multline*}
z=5625a^4c-1250a^3b^2-21600a^3c+4800a^2b^2-120000a^2bc-60000a^2c^2+\\
+24000ab^3+60000ab^2c+414720abc+1036800ac^2-10000b^4-\\
-81920b^3-38400b^2c-384000bc^2+160000c^3-2985984c^2.\end{multline*}
For a fixed $\gamma$ the curves $s(a,b,\gamma)=0$ and
$z(a,b,\gamma)=0$ divide the main triangle $\Delta_\gamma$ into
domains. In contrast to the behavior of the curve
$s(a,b,\gamma)=0$ the behavior of the curve $z(a,b,\gamma)=0$
inside $\Delta_\gamma$ vary with the change of $\gamma$. We will
describe the partition of $\Delta_\gamma$ into domains by curves
$s(a,b,\gamma)=0$ and $z(a,b,\gamma)=0$ and bifurcations of this
partition. Also we will indicate for each domain the passport of
the corresponding snake. \pmn Below are enumerated passports of
all 16 snake of order 5.
$$\begin{array}{rrrrrrr}
1)\,1,3,2,5,4&&2)\,1,4,2,5,3&&3)\,1,4,3,5,2&&
4)\,1,5,2,4,3\\ 5)\,1,5,3,4,2&&6)\,2,3,1,5,4&&7)\,2,4,1,5,3&&8)\,2,4,3,5,1\\
9)\,2,5,1,4,3&&10)\,2,5,3,4,1&&11)\,3,4,1,5,2&&12)\,3,4,2,5,1\\
13)\,3,5,1,4,2&&14)\,3,5,2,4,1&&15)\,4,5,1,3,2&&16)\,4,5,2,3,1\end{array}$$
In what follows we will schematically demonstrate the partition of
the main triangle into domains. Curves $z(a,b,\gamma)=0$ will be
dotted.

\subsection{$1>c>512/625\approx 0.8192$} If $\gamma$ belongs to this
interval, then the curve $z(a,b,\gamma)=0$ doesn't intersect the
main triangle. A point in a domain defines the Morse polynomial.
The plot of this polynomial is a snake and the number of its
passport is indicated inside the domain.
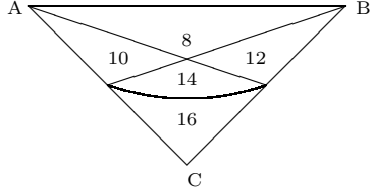
\[\begin{picture}(140,80) \put(10,70){\line(1,0){120}}
\put(10,70){\line(1,-1){60}} \put(70,10){\line(1,1){60}}
\put(10,70){\line(3,-1){90}} \put(40,40){\line(3,1){90}}
\qbezier(40,40)(70,30)(100,40) \put(2,67){\scriptsize A}
\put(134,67){\scriptsize B} \put(70,2){\scriptsize C}
\put(40,48){\scriptsize 10} \put(68,55){\scriptsize 8}
\put(92,48){\scriptsize 12} \put(66,40){\scriptsize 14}
\put(66,25){\scriptsize 16} \end{picture}\] If $\gamma=512/625$,
then the curve $z(a,b,\gamma)=0$ passes through the point $A$ --- the
vertex of the main triangle.

\subsection{$512/625>c>432/625\approx 0.6912$}
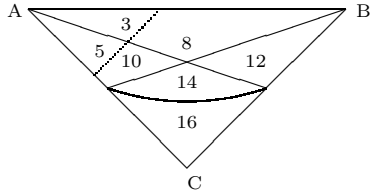
\[\begin{picture}(140,80) \put(10,70){\line(1,0){120}}
\put(10,70){\line(1,-1){60}} \put(70,10){\line(1,1){60}}
\put(10,70){\line(3,-1){90}} \put(40,40){\line(3,1){90}}
\qbezier(40,40)(70,30)(100,40) \qbezier[25](35,45)(45,55)(60,70)
\put(2,67){\scriptsize A} \put(134,67){\scriptsize B}
\put(70,2){\scriptsize C} \put(45,48){\scriptsize 10}
\put(68,55){\scriptsize 8} \put(92,48){\scriptsize 12}
\put(66,40){\scriptsize 14} \put(66,25){\scriptsize 16}
\put(45,62){\scriptsize 3} \put(35,52){\scriptsize 5}
\end{picture}\] If $\gamma=432/625$, then the curve $z(a,b,\gamma)=0$
passes through the point $B$  --- the vertex of the main triangle.

\subsection{$432/625>c>52488/72125\approx 0.6718$}
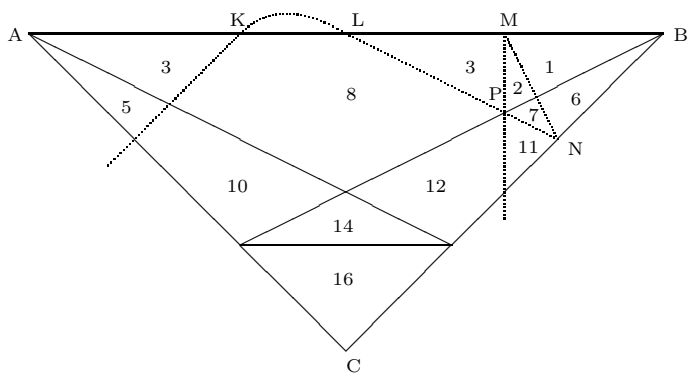
\[\begin{picture}(260,140) \put(10,130){\line(1,0){240}}
\put(10,130){\line(1,-1){120}} \put(10,130){\line(2,-1){160}}
\put(90,50){\line(2,1){160}} \put(90,50){\line(1,0){80}}
\put(130,10){\line(1,1){120}}
\qbezier[35](190,130)(200,110)(210,90)
\qbezier[50](190,130)(190,95)(190,60)
\qbezier[60](210,90)(170,110)(130,130)
\qbezier[50](90,130)(70,110)(40,80)
\qbezier[30](90,130)(105,145)(130,130) \put(2,127){\scriptsize A}
\put(254,127){\scriptsize B} \put(130,2){\scriptsize C}
\put(86,133){\scriptsize K} \put(132,133){\scriptsize L}
\put(188,133){\scriptsize M} \put(214,84){\scriptsize N}
\put(184,105){\scriptsize P} \put(60,115){\scriptsize 3}
\put(130,105){\scriptsize 8}  \put(175,115){\scriptsize 3}
\put(205,115){\scriptsize 1} \put(193,107){\scriptsize 2}
\put(45,100){\scriptsize 5} \put(85,70){\scriptsize 10}
\put(125,55){\scriptsize 14} \put(125,35){\scriptsize 16}
\put(160,70){\scriptsize 12} \put(195,85){\scriptsize 11}
\put(199,97){\scriptsize 7} \put(215,103){\scriptsize 6}
\end{picture}\]  Here points $M$ and $N$ lie on the sides of the
main triangle and both are cusps of the curve $z(a,b,\gamma)=0$.
The point $P$ lies on the curve $s(a,b,\gamma)=0$ and is an
self-intersection point of the curve $z$. The arc $KL$ disappears
for $\gamma=52488/72125$, i.e. the curve $z(a,b,\gamma)=0$ is
tangent (from inside) to the side $AB$ of the main triangle. With
the decrease of $\gamma$ two domains with numbers 3 merge into
one.

\subsection{$52488/72125>c>\approx 0.57613$}
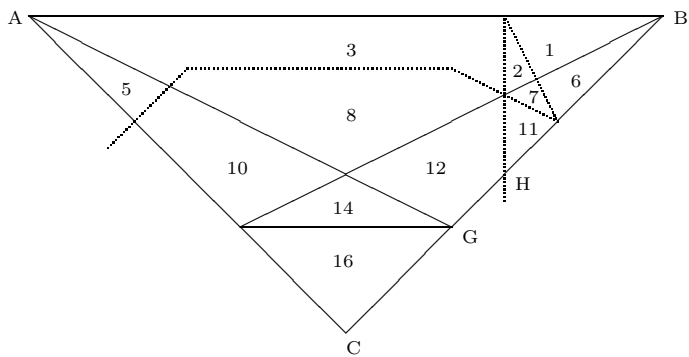
\[\begin{picture}(260,140) \put(10,130){\line(1,0){240}}
\put(10,130){\line(1,-1){120}} \put(10,130){\line(2,-1){160}}
\put(90,50){\line(2,1){160}} \put(90,50){\line(1,0){80}}
\put(130,10){\line(1,1){120}}
\qbezier[35](190,130)(200,110)(210,90)
\qbezier[50](190,130)(190,95)(190,60)
\qbezier[30](210,90)(190,100)(170,110)
\qbezier[30](70,110)(55,95)(40,80)
\qbezier[60](70,110)(120,110)(170,110)\put(2,127){\scriptsize A}
\put(254,127){\scriptsize B} \put(130,2){\scriptsize C}
\put(194,64){\scriptsize H} \put(174,44){\scriptsize G}
\put(130,115){\scriptsize 3} \put(130,90){\scriptsize 8}
\put(205,115){\scriptsize 1} \put(193,107){\scriptsize 2}
\put(45,100){\scriptsize 5} \put(85,70){\scriptsize 10}
\put(125,55){\scriptsize 14} \put(125,35){\scriptsize 16}
\put(160,70){\scriptsize 12} \put(195,85){\scriptsize 11}
\put(199,97){\scriptsize 7} \put(215,103){\scriptsize 6}
\end{picture}\] Points $G$ and $H$ became one for
$\gamma\approx 0.57613$. This value of $\gamma$ is a root
of a polynomial of degree 6, where the ratio of lowest coefficient to the
highest is $2^{52}/5^{15}$.

\subsection{$0.57613>c>1024/1875\approx 0.5461$}
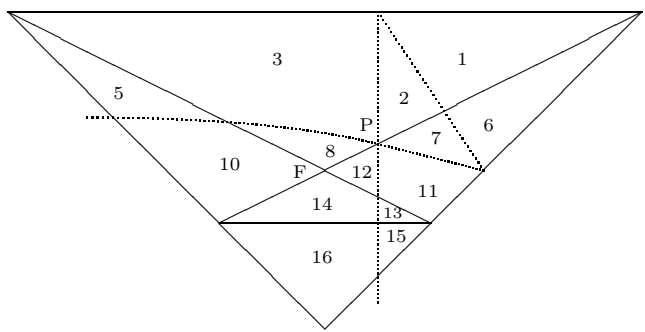
\[\begin{picture}(260,140) \put(10,130){\line(1,0){240}}
\put(10,130){\line(1,-1){120}} \put(10,130){\line(2,-1){160}}
\put(90,50){\line(2,1){160}} \put(90,50){\line(1,0){80}}
\put(130,10){\line(1,1){120}} \qbezier[70](40,90)(110,90)(150,80)
\qbezier[30](150,80)(170,75)(190,70)
\qbezier[40](190,70)(170,100)(150,130)
\qbezier[70](150,130)(150,75)(150,20) \put(118,67){\scriptsize F}
\put(143,85){\scriptsize P} \put(50,97){\scriptsize 5}
\put(110,110){\scriptsize 3} \put(158,95){\scriptsize 2}
\put(180,110){\scriptsize 1} \put(90,70){\scriptsize 10}
\put(130,75){\scriptsize 8} \put(140,67){\scriptsize 12}
\put(165,60){\scriptsize 11} \put(170,80){\scriptsize 7}
\put(190,85){\scriptsize 6} \put(125,55){\scriptsize 14}
\put(125,35){\scriptsize 16} \put(153,43){\scriptsize 15}
\put(152,52){\tiny 13}  \end{picture}\]  When $c=1024/1875$, then
points $F$ and $P$ merge and domains 8 and 12 disappear.

\subsection{$0<c<1024/1875$} Here configuration of domains in the main
triangle does not change:
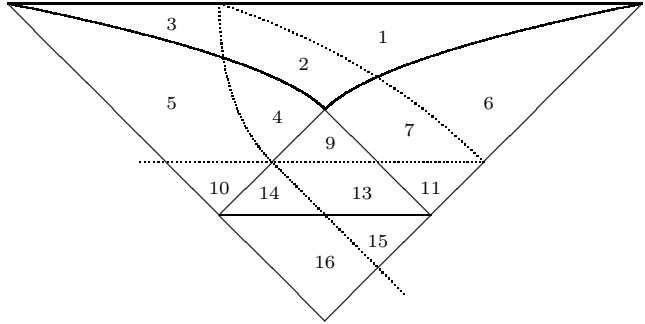
\[\begin{picture}(260,140) \put(10,130){\line(1,0){240}}
\put(10,130){\line(1,-1){120}} \put(130,10){\line(1,1){120}}
\put(90,50){\line(1,0){80}} \put(90,50){\line(1,1){40}}
\put(130,90){\line(1,-1){40}} \qbezier(130,90)(150,110)(250,130)
\qbezier(130,90)(110,110)(10,130)
\qbezier[80](60,70)(130,70)(190,70)
\qbezier[70](190,70)(150,110)(90,130)
\qbezier[50](90,130)(90,90)(110,70)
\qbezier[50](110,70)(130,50)(160,20) \put(70,120){\scriptsize 3}
\put(70,90){\scriptsize 5} \put(120,105){\scriptsize 2}
\put(150,115){\scriptsize 1} \put(190,90){\scriptsize 6}
\put(160,80){\scriptsize 7} \put(110,85){\scriptsize 4}
\put(86,58){\scriptsize 10} \put(105,56){\scriptsize 14}
\put(140,56){\scriptsize 13} \put(166,58){\scriptsize 11}
\put(146,38){\scriptsize 15} \put(126,30){\scriptsize 16}
\put(130,75){\scriptsize 9}
\end{picture}\]

\section{Supplement: the construction of Morse polynomials}
\pn Let $0=x_0<x_1<\ldots<x_k$ be critical points of a real Morse
polynomial $p$. Then $p'=x(x-x_1)\ldots(x-x_k)$ and
$$p=\int_0^x p'(t)\,dt.$$ The plot of $p'$ is of the form (for odd $k$):
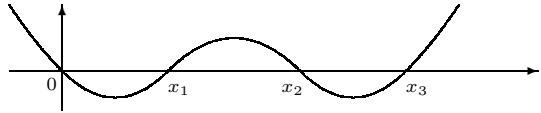
\[\begin{picture}(170,50) \put(0,20){\vector(1,0){200}}
\put(20,5){\vector(0,1){40}} \qbezier(20,20)(40,0)(60,20)
\qbezier(60,20)(85,45)(110,20) \qbezier(110,20)(130,0)(150,20)
\qbezier(150,20)(160,30)(170,45) \qbezier(20,20)(10,30)(0,45)
\put(14,13){\scriptsize 0} \put(60,12){\scriptsize $x_1$}
\put(103,12){\scriptsize $x_2$} \put(150,12){\scriptsize $x_3$}
\end{picture}\] We see the sequence of domains, bounded by the plot of
$p'$ and axis $OX$. Areas of these domains are differences of
consecutive critical values. Thus, the diminishing of
$x_{i+1}-x_i$ implies the diminishing of $|p(x_{i+1}-p(x_i)|$.
\begin{ex} Let us construct a Morse polynomial, whose plot has the
passport $(3,1,4,2)$. It will take several steps:
\begin{itemize}
\item roots $(0,1,2,3)$ --- passport $(4,2,3,1)$;
\item roots $(0,1,3,4)$ --- passport $(2,1,4,3)$;
\item roots $(0,1,3,5)$ --- passport $(3,2,4,1)$;
\item roots $(0,1,3,4.4)$ --- passport $(3,1,4,2)$.
\end{itemize} Analogously we can construct a 7-degree Morse polynomial
with the passport $(4,1,5,3,6,2)$. Its critical points are
$0,1,3,5,7,8.4$.
\end{ex}

\vspace{5mm}
\end{document}